\def\DATE{\relax}
\magnification=1100
\baselineskip=12.72pt
\voffset=.75in
\hoffset=.9in
\hsize=4.1in
\newdimen\hsizeGlobal
\hsizeGlobal=4.1in%
\vsize=7.05in
\parindent=.166666in
\pretolerance=500 \tolerance=1000 \brokenpenalty=5000

\footline={\vbox{\hsize=\hsizeGlobal\hfill{\rm\the\pageno}\hfill\llap{\sevenrm\DATE}}\hss}

\def\note#1{%
  \hfuzz=50pt%
  \vadjust{%
    \setbox1=\vtop{%
      \hsize 3cm\parindent=0pt\eightpoints\baselineskip=9pt%
      \rightskip=4mm plus 4mm\raggedright#1\hss%
      }%
    \hbox{\kern-4cm\smash{\box1}\hss\par}%
    }%
  \hfuzz=0pt
  }
\def\note#1{\relax}

\def\anote#1#2#3{\smash{\kern#1in{\raise#2in\hbox{#3}}}%
  \nointerlineskip}     % permet de mettre des anotations
                        % dans les graphiques inclus avec
                        % \special

\newcount\equanumber
\equanumber=0
\newcount\sectionnumber
\sectionnumber=0
\newcount\subsectionnumber
\subsectionnumber=0
\newcount\snumber  %statement number
\snumber=0

\def\section#1{%
  \subsectionnumber=0%
  \snumber=0%
  \equanumber=0%
  \advance\sectionnumber by 1%
  \noindent{\bf \the\sectionnumber .~#1.~}%
}%
\def\subsection#1{%
  \advance\subsectionnumber by 1%
  \snumber=0%
  \equanumber=0%
  \noindent{\bf \the\sectionnumber .\the\subsectionnumber .~#1.~}%
}%
\def\prevs{\the\sectionnumber .\the\subsectionnumber .\the\snumber }%previous statement
\long\def\Definition#1{%
  \global\advance\snumber by 1%
  \bigskip
  \noindent{\bf Definition~\prevs .}%
  \quad{\it#1}%
}

\long\def\Corollary#1{%
  \global\advance\snumber by 1%
  \bigskip
  \noindent{\bf Corollary~\prevs .}%
  \quad{\it#1}%
}%
\long\def\Lemma#1{%
  \global\advance\snumber by 1%
  \bigskip
  \noindent{\bf Lemma~\prevs .}%
  \quad{\it#1}%
}%
\def\Proof{\noindent{\bf Proof.~}}
\long\def\Proposition#1{%
  \advance\snumber by 1%
  \bigskip
  \noindent{\bf Proposition~\prevs .}%
  \quad{\it#1}%
}%
\long\def\Remark#1{%
  \bigskip
  \noindent{\bf Remark.~}#1%
}%
\long\def\remark#1{%
  \advance\snumber by1%
  \bigskip
  \noindent{\bf Remark~\prevs .}\quad#1%
}%
\long\def\Theorem#1{%
  \advance\snumber by 1%
  \bigskip
  \noindent{\bf Theorem~\prevs .}%
  \quad{\it#1}%
}%
\long\def\Statement#1{%
  \advance\snumber by 1%
  \bigskip
  \noindent{\bf Statement~\prevs .}%
  \quad{\it#1}%
}%
% labelling statements
\def\ifundefined#1{\expandafter\ifx\csname#1\endcsname\relax}
\def\labeldef#1{\global\expandafter\edef\csname#1\endcsname{\prevs}}
\def\labelref#1{\expandafter\csname#1\endcsname}
\def\label#1{\ifundefined{#1}\labeldef{#1}\note{$<$#1$>$}\else\labelref{#1}\fi}

% labeling and numbering of the equations
\def\preveq{(\the\sectionnumber .\the\subsectionnumber .\the\equanumber)}
\def\neq{\global\advance\equanumber by 1\eqno{\preveq}}

\def\ifundefined#1{\expandafter\ifx\csname#1\endcsname\relax}

%version of equadef with marginal note
%\def\equadef#1{\global\advance\equanumber by 1%
%  \global\expandafter\edef\csname#1\endcsname{\preveq}%
%  \setbox1=\hbox{\rm\hskip .1in[#1]}\dp1=0pt\ht1=0pt\wd1=0pt%
%  \preveq\box1}
%version of equadef with no marginal note
\def\equadef#1{\global\advance\equanumber by 1%
  \global\expandafter\edef\csname#1\endcsname{\preveq}%
  \preveq}

\def\equaref#1{\expandafter\csname#1\endcsname}

\def\equa#1{%
  \ifundefined{#1}%
    \equadef{#1}%
  \else\equaref{#1}\fi}

\font\eightrm=cmr8%
\font\sixrm=cmr6%

\font\eightsl=cmsl8%

\font\eightbf=cmb8%

\font\eighti=cmmi8%
\font\sixi=cmmi6%

\font\eightsy=cmsy8%
\font\sixsy=cmsy6%

\font\eightex=cmex8%
\font\sixex=cmex6%
\font\fiveex=cmex5%

\font\eightit=cmti8%

\font\eighttt=cmtt8%

\font\tenbb=msbm10%
\font\eightbb=msbm8%
\font\sevenbb=msbm7%
\font\sixbb=msbm6%
\font\fivebb=msbm5%
\newfam\bbfam  \textfont\bbfam=\tenbb  \scriptfont\bbfam=\sevenbb  \scriptscriptfont\bbfam=\fivebb%

\font\tenbbm=bbm10

\font\tencmssi=cmssi10%
\font\sevencmssi=cmssi7%
\font\fivecmssi=cmssi5%
\newfam\ssfam  \textfont\ssfam=\tencmssi  \scriptfont\ssfam=\sevencmssi  \scriptscriptfont\ssfam=\fivecmssi%

\font\tenfrak=cmfrak10%
\font\eightfrak=cmfrak8%
\font\sevenfrak=cmfrak7%
\font\sixfrak=cmfrak6%
\font\fivefrak=cmfrak5%
\newfam\frakfam  \textfont\frakfam=\tenfrak  \scriptfont\frakfam=\sevenfrak  \scriptscriptfont\frakfam=\fivefrak%
\def\frak{\fam\frakfam\tenfrak}%

% ams pour les signes \leq et \geq
\font\tenmsam=msam10%
\font\eightmsam=msam8%
\font\sevenmsam=msam7%
\font\sixmsam=msam6%
\font\fivemsam=msam5%

\def\bb{\fam\bbfam\tenbb}%

\def\hexdigit#1{\ifnum#1<10 \number#1\else%
  \ifnum#1=10 A\else\ifnum#1=11 B\else\ifnum#1=12 C\else%
  \ifnum#1=13 D\else\ifnum#1=14 E\else\ifnum#1=15 F\fi%
  \fi\fi\fi\fi\fi\fi}
\newfam\msamfam  \textfont\msamfam=\tenmsam  \scriptfont\msamfam=\sevenmsam  \scriptscriptfont\msamfam=\fivemsam%
\def\msam{\msamfam\tenmsam}%
\mathchardef\leq"3\hexdigit\msamfam 36%
\mathchardef\geq"3\hexdigit\msamfam 3E%

\font\tentt=cmtt11%
\font\seventt=cmtt9%
\textfont\ttfam=\tentt
\scriptfont7=\seventt%
\def\tt{\fam\ttfam\tentt}%

\def\eightpoints{%
\def\rm{\fam0\eightrm}%
\textfont0=\eightrm   \scriptfont0=\sixrm   \scriptscriptfont0=\fiverm%
\textfont1=\eighti    \scriptfont1=\sixi    \scriptscriptfont1=\fivei%
\textfont2=\eightsy   \scriptfont2=\sixsy   \scriptscriptfont2=\fivesy%
\textfont3=\eightex   \scriptfont3=\sixex   \scriptscriptfont3=\fiveex%
\textfont\itfam=\eightit  \def\it{\fam\itfam\eightit}%
\textfont\slfam=\eightsl  \def\sl{\fam\slfam\eightsl}%
\textfont\ttfam=\eighttt  \def\tt{\fam\ttfam\eighttt}%
\textfont\bffam=\eightbf  \def\bf{\fam\bffam\eightbf}%

\textfont\frakfam=\eightfrak  \scriptfont\frakfam=\sixfrak \scriptscriptfont\frakfam=\fivefrak  \def\frak{\fam\frakfam\eightfrak}%
\textfont\bbfam=\eightbb      \scriptfont\bbfam=\sixbb     \scriptscriptfont\bbfam=\fivebb      \def\bb{\fam\bbfam\eightbb}%
\textfont\msamfam=\eightmsam  \scriptfont\msamfam=\sixmsam \scriptscriptfont\msamfam=\fivemsam  \def\msam{\msamfam\eightmsam}

\rm%
}

\mathchardef\lsim"3\hexdigit\msamfam 2E%
\mathchardef\gsim"3\hexdigit\msamfam 26%

\def\d{\,{\rm d}}

\def\ds{\displaystyle}
\long\def\DoNotPrint#1{\relax}
\def\<{\langle}
\def\>{\rangle}
\def\d{\,{\rm d}}
\def\ds{\displaystyle}
\def\eqd{\,{\buildrel{\rm d}\over =}\,}
\def\fixedRef#1{#1\note{fixedRef}}

\def\prob#1{{\rm P}\{\,#1\,\}}
\def\Prob#1{{\rm P}\bigl\{\,#1\,\bigr\}}
\def\probb#1{{\rm P}\Bigl\{\,#1\,\Bigr\}}
\def\qed{{\vrule height .9ex width .8ex depth -.1ex}}

\def\oH{\overline H}
\def\One{\hbox{\tenbbm 1}}

\def\boc{\note{{\bf BoC}\hskip-11pt\setbox1=\hbox{$\Bigg\downarrow$}%
         \dp1=0pt\ht1=0pt\ht1=0pt\leavevmode\raise -20pt\box1}}
\def\eoc{\note{{\bf EoC}\hskip-11pt\setbox1=\hbox{$\Bigg\uparrow$}%
         \dp1=0pt\ht1=0pt\ht1=0pt\leavevmode\raise 20pt\box1}}

\def\One{\hbox{\tenbbm 1}}

\def\calR{{\cal R}}
\def\calS{{\cal S}}
\def\calT{{\cal T}}

\def\E{\hbox{\rm E}}

\def\R{{\cal R}}
\def\T{{\cal T}}

\def\MM{{\bb M}\kern .4pt}
\def\NN{{\bb N}\kern .5pt}
\def\RR{{\bb R}}

\pageno=1

\centerline{\bf A CONDITIONAL LIMIT THEOREM}
\centerline{\bf FOR A BIVARIATE REPRESENTATION}
\centerline{\bf OF A UNIVARIATE RANDOM VARIABLE}
\centerline{\bf AND CONDITIONAL EXTREME VALUES}

\bigskip
 
\centerline{Ph.\ Barbe$^{(1)}$ and Miriam Isabel Seifert$^{(2)}$}
\centerline{${}^{(1)}$CNRS {\sevenrm(UMR {\eightrm 8088})}, ${}^{(2)}$Helmut Schmidt Universit\"at}

%\bigskip
 
{\narrower
\baselineskip=9pt\parindent=0pt\eightpoints

\bigskip

{\bf Abstract.} We consider a real random variable $X$ represented through a
random pair $(R,T)$ in $\RR^2$ and a deterministic function $u$ as 
$X=Ru(T)$. Under 
some additional assumptions, we prove a limit theorem for $(R,T)$ given $X>x$,
as $x$ tends to infinity. As a consequence, we derive conditional limit 
theorems for random pairs $(X,Y)=(Ru(T),Rv(T))$ given that $X$ is large. 
These results imply earlier ones which were obtained in the literature
under stronger assumptions.

\bigskip


\bigskip
 

}

\bigskip
%%%%%%%%%%%%%%%%%%%%%%%%%%%%%%%%%%%%%%%%%%%%%%%%%%%%%%%%%%%%%%%%
%%%%%%%%%%%%%%%%%%%%%%%%%%%%%%%%%%%%%%%%%%%%%%%%%%%%%%%%%%%%%%%%

%%%%%%%%%%%%%%%%%%%%%%%%%%%%%%%%%%%%%%%%%%%%%%%%%%%%%%%%%%%%%%%%
\def\prevs{\the\sectionnumber.\the\snumber }%previous statement
\def\preveq{(\the\sectionnumber.\the\equanumber)}

\section{Introduction}
The purpose of this paper is to clarify some conditional limit theorems on 
bivariate vectors given that one of the component is large. The significance
of such limit theorems stems from their applications in multivariate extreme
value theory, where one is interested in both making statistical inference on a
system given that a component has an extreme behavior and understanding
the dependence structure between extreme events. These conditional theorems 
provide the theoretical support in the study of extremal behavior of random
vectors in the conditional extreme value models introduced
by Heffernan and Tawn (2004) and Heffernan and Resnick (2007), Das and Resnick
(2011), as well as for studying estimators in statistical applications as
done by Foug\`eres and Soulier (2012).

Following these authors we are interested in a generalization of elliptically
distributed random vectors, namely, random vectors $(X,Y)$ with
representation $\bigl(Ru(T),Rv(T)\bigr)$, where $u$ and $v$ are 
deterministic functions, $R$ and $T$ are independent real random variables, 
and the distribution of $R$ is in the Gumbel 
max-domain of attraction (see Berman, 1983; Foug\`eres and
Soulier, 2010; Hashorva, 2012; Seifert, 2012). For elliptical random variables,
$R$ is the radial component and $T$ the angular distribution. However, in our
more general setting, the map $(R,T)\mapsto (X,Y)$ may not be one-to-one.

Beyond immediate applications to extreme value theory, our results have 
bearing to the description of the convex hull of samples and related 
problems which are in part driven by extreme value theory.

The novelty of our paper is to show that a conditional limit theorem 
for $(X,Y)$ given that $X$ is large is not intrinsically about the pair $(X,Y)$
but about the representation of the single variable $X$ in terms of the pair
$(R,T)$. This approach allows us to recover previous results under minimal 
assumptions, to provide a better understanding of the earlier work, and, 
through more versatile assumptions, to widen the applicability of this
model.

Throughout the paper, $R$ is a real random variable, so that
$R\bigl(u(T),v(T)\bigr)$ means $\bigl(Ru(T),Rv(T)\bigr)$.

\bigskip

%%%%%%%%%%%%%%%%%%%%%%%%%%%%%%%%%%%%%%%%%%%%%%%%%%%%%%%%%%%%%%%%%%%%%%
%
%%%%%%%%%%%%%%%%%%%%%%%%%%%%%%%%%%%%%%%%%%%%%%%%%%%%%%%%%%%%%%%%%%%%%%
\def\prevs{\the\sectionnumber .\the\snumber }%previous statement
\def\preveq{(\the\sectionnumber .\the\equanumber)}

\section{Main result}
In this section, we are interested in random 
variables $X$ which are represented as $X=Ru(T)$,
and conditional limit theorems for properly normalized $(R,T)$ given $X>x$ 
as $x$ tends to infinity. In the next section, equipped with such a 
conditional limit theorem we will use some
continuous mapping argument to derive a conditional limit theorem 
for properly normalized $Y=Rv(T)$ given $X>x$ as $x$ tends to infinity.

We write $H$ for the cumulative distribution function of $R$, and $\oH$ for
the survival function $1-H$. We assume that $T$ has a density $g(t)$. 

We will use the following assumptions.

\bigskip

\noindent{\bf Assumption 1.} {\it\ %
  The survival function $\oH$ of $R$ is in the class $\Gamma(\psi)$,}

\bigskip

\noindent
meaning that there exists an ultimately positive function $\psi$ such that
for any fixed real number $\lambda$,
$$
  \lim_{x\to\infty} {\oH\bigl(x+\psi(x)\lambda\bigr)\over\oH(x)} 
  = e^{-\lambda} \, .
$$
This property is equivalent to $H$ belonging to the max-domain of 
attraction of the Gumbel distribution
(de Haan, 1970; Resnick, 2007). The function $\psi$ is unique up to
asymptotic equivalence, and, necessarily, $\psi(x)=o(x)$ at infinity.

\bigskip

\noindent{\bf Assumption 2.} {\it\ %
  There exists a $t_0$ such that $u(t_0)=1$ and for any $\epsilon$ positive, 
  $\sup_{t-t_0>\epsilon} u(t)<1$. Moreover, the function
  $$
    \tilde u(s)=u(t_0)-u(t_0+s)
  $$
  is regularly varying at $0+$ with positive index $\kappa$,
}

\bigskip

\noindent meaning that for any positive $\lambda$,
$$
  \lim_{s\to 0+} {\tilde u(\lambda s)\over \tilde u(s)} = \lambda^\kappa \, .
$$
The first part of assumption 2 asserts that on the right of $t_0$, the
function $u$ has a unique maximum at $t_0$ 
and that for 
$u(t)$ to be close to $1$, we must have $t$ close to $t_0$.

Since $\psi(x)=o(x)$ at infinity and $\tilde u$ is regularly varying with
positive index, there exists an ultimately positive function $\phi$ such that
$$
  \tilde u\circ \phi(x)\sim {\psi(x)\over x}
  \eqno{\equa{phiDef}}
$$
as $x$ tends to infinity, and $\lim_{x\to \infty}\phi(x)=0$. 

We will also use the notation
$$ 
  \tilde g(s)=g(t_0+s) \, ,
$$
and assume that

\bigskip

\noindent {\bf Assumption 3.}{\it\ \ %
   The density $\tilde g$ of $T-t_0$ is regularly varying at $0+$
   with index $\tau>-1$.
}

\bigskip

Since $\tilde g$ is locally integrable, its index of regular variation must be
at least $-1$. Furthermore, if $\tilde g$ is positive and continuous at $0$, 
then $\tau$ vanishes.

To keep track of the notation, note that whenever a function has a tilde, 
it means that it is regularly varying at $0$.

Our main result is the following conditional limit theorem for $(R,T)$ given
$X>x$ and $T>t_0$, as $x$ tends to infinity. We will see in the next section
how the conditioning by $T>t_0$ may be removed under additional assumptions.

\bigskip

\Theorem{\label{mainTh}%
  Let $X=Ru(T)$. Under assumptions 1, 2 and 3, the conditional distribution of
  $$
    \Bigl( {R-x\over \psi(x)},{T-t_0\over \phi(x)}\Bigr)
  $$
  given $X>x$ and $T>t_0$ converges weakly$*$, as $x$ tends to infinity, 
  to the measure whose density with respect to the Lebesgue measure is
  $$
    {\kappa\over \Gamma\Bigl({\ds 1+\tau\over\ds\kappa}\Bigr) }
    t^\tau e^{-r}
    \One\{\, 0<t<r^{1/\kappa} \,\}
    \eqno{\equa{limitDensity}}
  $$
  as $x$ tends to infinity. Furthermore,
  $$
    \prob{X>x \, ;\, T>t_0 }
    \sim \phi(x)\tilde g\circ\phi(x) \overline H(x) 
       {1\over\kappa}\Gamma\Bigl({1+\tau\over \kappa}\Bigr)
  $$
  as $x$ tends to infinity.
}

\bigskip

Some heuristic arguments explaining why Theorem \mainTh\ may be true are given
at the beginning of section 5.

As a function defined on some right neighborhood of the origin, $\tilde u$ has
an asymptotic inverse $\tilde u^\leftarrow$ such 
that $\tilde u\circ\tilde u^\leftarrow(s)\sim s$ as $s$ tends to $0$ (see
Bingham, Goldie and Teugels, 1989, \S 1.5.7). Thus, 
$\phi(x)\sim \tilde u^\leftarrow \bigl( \psi(x)/x\bigr)$ and
$$
  (\phi g\circ \phi)(x)
  \sim (\tilde u^\leftarrow \tilde g\circ \tilde u^\leftarrow)(x)
$$
as $x$ tends to infinity. Therefore, we may view $\tilde \phi\tilde g\circ\phi$
as a function of $\psi(x)/x$ which is then regularly varying of 
index $(1+\tau)/\kappa$ in terms of the argument $\psi(x)/x$.

\bigskip

\section{Two-sided extensions}
In some applications it is desirable to have analogues of Theorem \mainTh\
when the conditioning involves only the event $X>x$. Under
two-sided conditions on the behavior of $\tilde u$ and $\tilde g$ near $t_0$, 
such extensions present no conceptual difficulty. To illustrate this assertion,
we present two such extensions, relying on the following two-sided
versions of assumptions 2 and 3.

\bigskip

\noindent{\bf Assumption 4.} {\it\ \ %
  There exists a $t_0$ such that $u(t_0)=1$ and for any $\epsilon$ positive,
  $\sup_{|t-t_0|>\epsilon} u(t)<1$. Moreover, the function $\tilde u$ is regularly
  varying at $0-$ and $0+$ with respective positive indices $\kappa_-$ 
  and $\kappa_+$.
}

\bigskip

\noindent
The second part of assumption 4 signifies that for any given sign $\sigma$ in 
$\{\, -,+\,\}$ and any positive $\lambda$
$$
  \lim_{s\to 0+} {\tilde u(\sigma\lambda s)\over \tilde u(\sigma s)}
  =\lambda^{\kappa_\sigma} \, .
$$

Similarly, we strengthen assumption 3 as follows.

\bigskip

\noindent{\bf Assumption 5.} {\it\ \ %
  $\tilde g$ is regularly varying at $0-$ and $0+$ with respective indices
  $\tau_-$ and $\tau_+$, both these indices being greater than $-1$.
}

\bigskip

Equipped with these two-sided hypotheses, we define, as in \phiDef, for each
sign $\sigma$, an ultimately positive function $\phi_\sigma$ such that
$$
  \tilde u \bigl( \sigma\phi_\sigma(x)\bigr)\sim {\psi(x)\over x}
$$
as $x$ tends to infinity, and $\lim_{x\to\infty}\phi_\sigma(x)=0$. In order
to describe the contributions of both sides of $t_0$ to the asymptotic
behavior of $(R,T)$, we further suppose the following.

\bigskip

\noindent{\bf Assumption 6.}{\it\ \ %
  For any sign $\sigma$,
  $$
    p_\sigma=\lim_{x\to\infty} 
    {\phi_\sigma \tilde g( \sigma\phi_\sigma) \over
     \phi_- \tilde g( -\phi_-)+\phi_+ \tilde g(\phi_+) }(x)
  $$
  exists.
}

\bigskip

Both $p_-$ and $p_+$ are nonnegative
and their sum is $1$. They represent the contribution of the events $T<t_0$
and $T>t_0$ to the limiting conditional distribution of $T-t_0$ given $X>x$.
Considering $\phi_\sigma \tilde g(\sigma \phi_\sigma)$ as a regularly varying
function of $\psi(x)/x$ of index $(1+\tau_\sigma)/\kappa_\sigma$, we see that
if both $p_-$ and $p_+$ do not vanish, 
then $(1+\tau_+)/\kappa_+=(1+\tau_-)/\kappa_-$.

To state our results, we introduce the random sign
$$
  S={\rm sign} (T-t_0) \, .
$$
We consider also a random sign $\calS$ whose distribution is
$$
  \prob{ \calS=\sigma }
  = {\ds {\ds p_\sigma\over\ds \kappa_\sigma}
      \Gamma\Bigl({\ds1+\tau_\sigma\over\ds\kappa_\sigma}\Bigr) 
      \over\ds 
      {\ds p_-\over\ds \kappa_-}
      \Gamma\Bigl({\ds1+\tau_-\over\ds\kappa_-}\Bigr) 
      +{\ds p_+\over\ds \kappa_+}
      \Gamma\Bigl({\ds1+\tau_+\over\ds\kappa_+}\Bigr) 
    }
  \, ,\qquad \sigma\in\{\, -,+\,\} \, .
$$
Central to our two-sided extension is the following consequence of Theorem
\mainTh. This result is also of importance to understand how the results
in the next section, stated under one-sided assumptions and an extra
conditioning on $T>t_0$, can be extended with two-sided assumptions and
no conditioning on $T>t_0$.

\Proposition{\label{signConditional}
  Under assumptions 1, 4, 5 and 6, the conditional distribution of $S$
  given $X>x$ converges weakly$*$ to that of $\calS$.
}

\bigskip

\Proof The second assertion of Theorem \mainTh\ implies that for any 
sign $\sigma$,
$$
  \prob{X>x\,;\, S=\sigma}
  \sim \phi_\sigma (x)\tilde g\bigl(\sigma\phi_\sigma(x)\bigr)\oH(x) 
  {1\over\kappa_\sigma}
  \Gamma\Bigl( {1+\tau_\sigma\over\kappa_\sigma}\Bigr)
$$
as $x$ tends to infinity. The proposition then follows from the
formula
$$
  \prob{S=\sigma \mid X>x}
  = {\prob{X>x \,;\, S=\sigma}\over \prob{X>x\,;\, S=-} +\prob{X>x\,;\, S=+}}
  \, . 
  \eqno{\qed}
$$

\bigskip

We then define a random pair $(\calR,\calT_\calS)$ whose conditional 
distribution given $\calS=\sigma$ has density with respect to the 
Lebesgue measure
$$
  {\kappa_\sigma\over \Gamma\Bigl({\ds 1+\tau_\sigma\over\ds\kappa_\sigma}\Bigr)}
  t^{\tau_\sigma} e^{-r}\One\{\, 0<t<r^{1/\kappa_\sigma}\,\} \, .
$$

\Theorem{\label{thCondA}
  Under assumptions 1, 4, 5 and 6, the conditional distribution of 
  $$
    \Bigl( {R-x\over\psi(x)},{T-t_0\over\phi_S(x)}\Bigr)
  $$
  given $X>x$ converges weakly$*$ as $x$ tends to infinity to the
  distribution of $(\calR,\calS\calT_\calS)$.
}

\bigskip

The density of the limiting distribution can be written explicitly as
$$
  \sum_{\sigma\in\{-,+\}}
  |t|^\tau e^{-r} 
  { p_\sigma\One\{\, |t|^{\kappa_\sigma} <r\,:\, \sigma t>0\,\}
    \over 
    \ds {\ds p_-\over\ds \kappa_-}\Gamma\Bigl({\ds 1+\tau_-\over\ds\kappa_-}
                                        \Bigr)  
       +{\ds p_+\over\ds \kappa_+}\Gamma\Bigl({\ds 1+\tau_+\over\ds\kappa_+}
                                        \Bigr)
  } \, .
$$

\Proof For any Borel subset $A$ of $\RR^2$, we have
$$\displaylines{\qquad
  \probb{ \Bigl( {R-x\over\psi(x)},{T-t_0\over\phi_S(x)}\Bigr)\in A
  \Bigm| X>x}
  \hfill\cr\hfill
  {}=\sum_{\sigma\in\{-,+\}}   
  \probb{ \Bigl( {R-x\over\psi(x)},{T-t_0\over\phi_S(x)}\Bigr)\in A
  \Bigm| X>x\,;\, S=\sigma}
  \qquad\cr\hfill
  \prob{ S=\sigma \mid X>x} \, .
  \qquad\equa{thCondAa}\cr}
$$
Theorem \mainTh\ implies that the conditional distribution of
$$
  \Bigl( {R-x\over\psi(x)},\sigma {T-t_0\over\phi_\sigma(x)}\Bigr)
$$
given $X>x$ and $S=\sigma$ converges weakly$*$ to that of a random variable
$(\calR,\calT_\sigma)$ whose density with respect to the Lebesgue measure is
$$
  {\kappa_\sigma\over\Gamma\Bigl({\ds 1+\tau_\sigma\over\ds\kappa_\sigma}\Bigr)}
  t^{-\tau_\sigma}e^{-r} \One\{\, 0<t<r^{1/\kappa_\sigma}\,\} \, .
$$
Combining Proposition \signConditional\ and \thCondAa, we obtain that the
conditional distribution of
$$
  \Bigl( {R-x\over\psi(x)},{T-t_0\over \phi_S(x)}\Bigr)
$$
given $X>x$ converges weakly$*$ to that of $(\calR,\calS\calT_\calS)$.\hfill\qed

\bigskip

One may argue that the random norming of $T-t_0$ by $1/\phi_S(x)$ in Theorem
\thCondA\ would be better replaced by a deterministic one. This can be done,
defining 
$$
  \phi_*=\phi_++\phi_-
$$
and assuming

\bigskip

\noindent{\bf Assumption 7.}{\it\ \ 
  For any sign $\sigma$, the limit
  $\ds q_\sigma=\lim_{x\to\infty} {\ds\phi_\sigma\over\ds\phi_*}(x)$ exists.
}

\bigskip

We then have the following.

\Theorem{\label{thCondB}
  Under assumptions 1, 4--7, the conditional distribution of
  $$
    \Bigl( {R-x\over\psi(x)},{T-t_0\over \phi_*(x)}\Bigr)
  $$
  given $X>x$ converges weakly$*$ as $x$ tends to infinity to the distribution
  of $(\calR,q_\calS\calS\calT_\calS)$.
}

\bigskip

\noindent
Again, the limiting density can be made explicit if needed.

\bigskip

\Proof Given Proposition \signConditional\ and the definition of $p_\sigma$,
the conditional distribution of the random variable $\phi_\calS(x)/\phi_*(x)$
given $X>x$ converges weakly$*$ to that of $q_\calS$, and this convergence
holds jointly with the conditional convergence of
$$
  \Bigl( {R-x\over\psi(x)},{T-t_0\over\phi_S(x)}\Bigr) \, .
$$
The result follows.\hfill\qed

\bigskip

\section{Bivariate conditional limit theorems}
The purpose of this section is to use Theorem \mainTh\ to shed a new
light on previous results dealing with conditional bivariate distributions
given one extreme component.

To do so, we consider another random variable, $Y=Rv(T)$, under the 
conditional distribution given $X>x$ and $T>t_0$. Below, we will make precise
why we condition on both $X>x$ and $T>t_0$. However, the conditioning 
by $T>t_0$ can be easily removed by imposing the proper two-sided condition
and using the same arguments used to extend Theorem \mainTh\ to Theorems
\thCondA\ and \thCondB. In particular, removing the conditioning
by $T>t_0$ does not seem to add any insight on the problem. Thus, we choose
to keep this conditioning to keep the exposition concise.
We set
$$
  R_x= {R-x\over\psi(x)}
  \qquad\hbox{and}\qquad
  T_x={T-t_0\over\phi(x)} \, .
$$
Under the conditional distribution given $X>x$ and $T>t_0$, Theorem \mainTh\ 
asserts that $(R_x,T_x)$ converges in distribution to some $(\R,\T)$
whose density with respect to the Lebesgue measure is given by \limitDensity.

Similarly to $\tilde u$, define $\tilde v(s)=v(t_0)-v(t_0+s)$. Let 
us assume that

\bigskip

\noindent{\bf Assumption 8.} $\rho=v(t_0)$ is well defined and
  $\tilde v$ is regularly varying at $0+$, with nonnegative index $\delta$.

\bigskip

Note that $\delta=0$ is allowed; one could also look at what happens if
$\delta$ is negative, using the same technique but working directly with
$v(t_0+s)$ instead of $\tilde v(s)$; so the sign of $\delta$ does
not really matter, but we will take it nonnegative in order to see how
some known results follow from Theorem \mainTh.

We have
$$\eqalign{
  Y
  &{}= Rv(T)\cr
  &{}=\bigl(x+\psi(x)R_x\bigr) v\bigl(t_0+\phi(x)T_x\bigr) \cr
  &{}=\bigl(x+\psi(x)R_x\bigr) 
    \bigl(\rho-\tilde v\bigl(\phi(x)T_x\bigr)\bigr)\cr
  &{}=\rho x +\rho\psi(x)R_x - \bigl( x+\psi(x)R_x\bigr) \tilde v
    \bigl(\phi(x)T_x\bigr) \cr
  }
$$
Since $\psi(x)=o(x)$ and $\tilde v$ is regularly varying and both $R_x$
and $T_x$ remain bounded in probability, we obtain, when $T_x$ is nonnegative,
$$
  Y
  =\rho x+\rho \psi(x)R_x -T_x^\delta x\tilde v\circ\phi(x)
    \bigl(1+o(1)\bigr) \, .
  \eqno{\equa{YBasicA}}
$$

Using the Skorokhod-Dudley-Wichura theorem (see e.g.\ Dudley, 1989, 
sections 11.6 and 11.7), we can assume that
we have versions of $R_x$ and $T_x$ which converge almost surely to
$(\R,\T)$ on the events $\{\, X>x\,;\, T>t_0 \,\}$. We then obtain, under 
the conditional distribution given $X>x$ and $T>t_0$,
$$
  Y\eqd \rho x +\rho\psi(x)\R \bigl(1+o(1)\bigr) 
  - \T^\delta x\tilde v\circ\phi(x)
  \bigl(1+o(1)\bigr)
$$
as $x$ tends to infinity. Given \phiDef, this means
$$
  Y\eqd \rho x +\rho x \tilde u\circ \phi(x)\R \bigl(1+o(1)\bigr) 
  - \T^\delta x\tilde v\circ\phi(x)
  \bigl(1+o(1)\bigr) \, .
  \eqno{\equa{YBasic}}
$$

Recall that $\tilde u$ is regularly varying with index $\kappa$ and $\tilde v$
is regularly varying with index $\delta$, and that we have $\R>\T^\kappa$ 
almost surely. We can now vary the assumptions in several ways, which we state 
as examples.

\Remark We can now see what happens if we do not wish to condition on $T>t_0$.
We need to introduce the random sign $S={\rm sign}(T-t_0)$ and follow what
was done in section \fixedRef{3}. Identity \YBasicA\ becomes, with rather
obvious notation,
$$
  Y=\rho x+\rho \psi(x)R_x\bigl(1+o(1)\bigr)
  -|T_x|^{\delta_S}x \tilde v\bigl( S\phi_S(x)\bigr)\bigl(1+o(1)\bigr) \, .
$$
One then needs to discuss the behavior of $\tilde v$ on both sides ot $0$,
both in terms of regular variation and sign, and one can also discuss the 
possible replacement of $\phi_\calS$ by $\phi_*$. Such a discussion requires to 
distinguish very many cases and does not appear to bring further 
understanding. Thus we choose to state results that seems to be the most 
useful to specialize in applications.

\bigskip 

\noindent{\bf Example 1.} We assume that
$$
  \lim_{s\to 0+} \rho\tilde u(s)/\tilde v(s)=0 \, .
  \eqno{\equa{FougeresSoulier}}
$$
This is implied by Foug\`eres and Soulier's assumption that $\delta<\kappa$,
and it is also satisfied whenever $\rho$ is $0$.
Theorem \mainTh\ implies the following result which was proved under stronger
assumptions in Foug\`eres and Soulier (2010), up to the conditioning by $T>t_0$
which can be removed in using the same arguments as in the previous section.
Our proof shows that while this
result looks like a truly two-dimensional result, it is really two-dimensional
in $(R,T)$ but one-dimensional in $(X,Y)$.

\Corollary{\label{FS}
  Under the assumptions of Theorem \mainTh, assumption 8 and \FougeresSoulier, the 
  conditional distribution of
  $$
    \Bigl( {X-x\over \psi(x)},{Y-\rho x\over x\tilde v\circ\phi(x)}\Bigr)
  $$
  given $X>x$ and $T>t_0$ converges weakly$*$ to that 
  of $(\R-\T^\kappa,-\T^\delta)$ as $x$ tends to infinity.
}

\bigskip

\Proof \YBasic\ gives
$Y\eqd \rho x - \T^\delta x\tilde v\circ\phi(x)\bigl(1+o(1)\bigr)$,
and we have the convergence in distribution
$$
  {Y-\rho x\over x\tilde v\circ\phi(x)} \to -\T^\delta \, .
$$
This is the result.\hfill\qed

\bigskip

Corollary \FS\ makes it  quite clear why the function $H_{\eta,\tau}$
come up in Foug\`eres and Soulier (2010): this is what one gets from
Theorem \mainTh\ and the continuous mapping theorem, and it occurs
because of what the joint distribution of $(\R,\T)$ is.

\bigskip

\noindent{\bf Example 2.} Assume that
$$
  \lim_{s\to 0+} |\tilde u(s)/\tilde v(s)|=+\infty \, .
  \eqno{\equa{deltaGTKappaAssumption}}
$$
This is the case if $\delta>\kappa$ for instance.

\Corollary{\label{deltaGTkappa}
  Under the assumptions of Theorem \mainTh, assumption 8 
  and \deltaGTKappaAssumption, the conditional distribution of
  $$
    \Bigl( {X-x\over\psi(x)},{Y-\rho x\over \psi(x)}\Bigr)
  $$
  given $X>x$ and $T>t_0$ converges weakly$*$ to that 
  of $(\R-\T^\kappa,\rho\R)$ as $x$ tends to infinity.
}

\bigskip

\Proof It follows from \YBasic.\hfill\qed

\bigskip

Note that when $\rho$ vanishes, Corollary \deltaGTkappa\ yields a limiting 
distribution with degenerate second marginal. This means that in the 
conditional distribution $Y=o_P\bigl(\psi(x)\bigr)$ as $x$ tends to infinity.

\bigskip

\noindent{\bf Example 3.} Assume that
$$
  \lim_{s\to 0+} \tilde u(s)/\tilde v(s)=C\in \RR \, .
  \eqno{\equa{vuAsympAssumption}}
$$
When $C=0$, this is example 1.

\Corollary{\label{vuAsymp}
  Under the assumption of Theorem \mainTh, asumption 8 and \vuAsympAssumption, 
  the conditional distribution of
  $$
    \Bigl( {X-x\over\psi(x)},{Y-\rho x\over x\tilde v\circ \phi(x)}\Bigr)
  $$
  given $X>x$ and $T>t_0$ converges weakly$*$ to that 
  of $(\R-\T^\kappa,C \rho \R-\T^\delta)$.
}

\bigskip

\Proof It follows from \YBasic.\hfill\qed

\bigskip

If $C\not=0$, then Corollary \vuAsymp\ asserts as well that the conditional
distribution of 
$$
  \Bigl( {X-x\over\psi(x)},{Y-\rho x\over\psi(x)}\Bigr)
$$
given $X>x$ and $T>t_0$ converges weakly$*$ to that 
of $(\R-\T^\kappa,\rho\R-\T^\delta/C)$. This restatement gives example 2 at 
the limit when $C$ tends to infinity.

\bigskip

\noindent{\bf Example 4.} Assume now that
$$
  v(t)=(t-t_0+\rho)u(t)\ \hbox{ in a neighborhood of $t_0$.}
  \eqno{\equa{SeifertA}}
$$
Note that $\rho=v(t_0)$ as required in assumption 8.
In this case, $\tilde v(s)=(\rho+s)\tilde u(s)-su(t_0)$. This
identity shows that $\delta=\kappa\wedge 1$ if $\rho\not=0$ and $\delta=1$ if
$\rho=0$; therefore, we can be in any of
the cases covered by examples 1, 2 or 3: for instance, $\kappa>1$ or $\rho=0$
yield \FougeresSoulier; $\kappa<1$ and $\rho\not=0$ 
yield \vuAsympAssumption; and $\kappa=1$ and $\rho\not=0$ may yield any 
of \FougeresSoulier,\deltaGTKappaAssumption\ or \vuAsympAssumption. The 
question 
arises as to whether it is possible to have a unified normalization 
for $Y$ for its conditional distribution to converge.
The following result shows that with assumption \SeifertA, we cannot
anymore normalize $Y$ by some deterministic quantities independent of $\kappa$
and $\rho$. However, we can use a normalization which
involves $X$, as for instance Heffernan and Resnick (2007) did. Up to
the conditioning on $T>t_0$, the following
result was obtained by Seifert (2012) under stronger conditions.

\Corollary{\label{SeifertB}
  Under the assumptions of Theorem \mainTh, the conditional distribution of
  $$
    \Bigl( {X-x\over\psi(x)},{(Y/X)-\rho\over\phi(x)}\Bigr)
  $$
  given $X>x$ and $T>t_0$ converges weakly$*$ to that 
  of $(\R-\T^\kappa,\T)$ as $x$ tends to infinity.
}

\bigskip

\Proof We have
$$
  Y
  =R v(T)
  =X {v\over u}(T)
  =X(T-t_0+\rho) \, .
$$
Thus 
$$
  {(Y/X)-\rho\over \phi(x)}=T_x \, .
$$
The result follows since $T_x$ converges in distribution 
to $\T$ when $X>x$.\hfill\qed

\bigskip

In typical situations, $v$ is continuous and montone on a neighborhood
of $t_0$, while $u$ is continuous and monotone on a punctured neigborhood
of $t_0$ and $\kappa>\delta$, as assumed in Foug\`eres and Soulier (2010).
As shown in Seifert (2012), a suitable reparametrisation of $T$ yields
\SeifertA.

\bigskip

\noindent{\bf Example 5.} The previous example can be generalized in the
following way, yielding a somewhat exotic limiting behavior. Define the 
function $\theta(t)$ by the relation
$$
  v(t)=\theta(t)u(t) \, .
  \eqno{\equa{thetaDef}}
$$
and assume that for some nonnegative integer $n$, $\theta$ is $n$ times
differentiable and
$$
  \theta^{(j)}(t_0)=0 \hbox{ if } j=1,2,\ldots,n-1,\quad\hbox{ and }\quad
  \theta^{(n)}(t_0)\not= 0 \, . 
  \eqno{\equa{thetaAssumption}}
$$
Put differently, $n$ corresponds to the first nonvanish Taylor coefficient of
$(v/u)(t)-(v/u)(t_0)$.

\Corollary{\label{thetaCorollary}
  Under the assumptions of Theorem \mainTh\ and \thetaDef, \thetaAssumption,
  the conditional distribution of
  $$
    \Bigl( {X-x\over \psi(x)},{(Y/X)-\theta(t_0)\over \phi(x)^n}\Bigr)
  $$
  given $X>x$ and $T>t_0$ converges weakly$*$ to that 
  of $(\R-\T^\kappa,\T^n \theta^{(n)}(t_0)/n!)$ as $x$ tends to infinity.
}

\bigskip

\Proof Using \thetaDef,
$$
  Y/X=\theta(T)=\theta(t_0+\phi(x)T_x)
$$
as $x$ tends to infinity.
Since $\phi(x)$ tends to $0$ as $x$ tends to infinity, Taylor formulas and
the convergence in distribution of $T_x$ to $\T$ yield
$$
  Y/X\eqd\theta(t_0)+\phi(x)^n\T^n\theta^{(n)}(t_0)/n! \bigl(1+o(1)\bigr)
$$
as $x$ tends to infinity, which is the result.\hfill\qed

\bigskip

Of course, one could extend this example further in assuming that 
$\theta(t)-\theta(t_0)$
is regularly varying at $t_0$, and numerous other variations are possible.

\bigskip

To conclude, since all the results of this paper use basic regular 
variation theory, it is certain that a truly multivariate extension is 
possible. Such extension
is not unique for there exists various theories of multivariate regular
variation, beyond what is popular in extreme value theory; see for
instance the works of Mershaert and Scheffer (2001), the book by Vladimirov, 
Drozzinov and Zvialov (1988), and some pointers in
Bingham, Goldie and Teugels (1989). Which one is the most relevant seems
application dependent.

\bigskip

%%%%%%%%%%%%%%%%%%%%%%%%%%%%%%%%%%%%%%%%%%%%%%%%%%%%%%%%%%%%%%%%%%%%%%
% Proofs
%%%%%%%%%%%%%%%%%%%%%%%%%%%%%%%%%%%%%%%%%%%%%%%%%%%%%%%%%%%%%%%%%%%%%%

\section{Proof of Theorem \mainTh}
Before giving a formal proof, it is enlightening to
give an intuition on how this result was found and why it might
be true. We have $X=Ru(T)$. If $X>x$ and $x$ is large,
since $u$ is at most $1$ and $R$ has a light tail, we should expect $R$ to be
about $x$ and $u(T)$ about $1$, that is, $T$ about $t_0$; more 
precisely, since $\oH$ is in $\Gamma(\psi)$,
we should have, $R\approx x+\psi(x)\R$ for some $\R$ of order
$1$, and, hopefully, $T\approx t_0+\phi(x)\T$ for some function $\phi$
which tends to $0$ at infinity, and some $\T$ of order $1$. Moreover, if
$T>t_0$ then $\calT$ should be nonnegative. That would give
$$
  X=Ru(T)
  \approx\bigl(x+\psi(x)\R\bigr) u\bigl(t_0+\phi(x)\T\bigr)
  \eqno{\equa{XAsympt}}
$$
One should then look at $u$ near $t_0$, and so we define
$$
  \tilde u(s)=u(t_0)-u(t_0+s) \, .
$$
If this function is regularly varying at $0+$ with index $\kappa$, and 
since $u(t_0)=1$, we expect 
$$
  u\bigl(t_0+\phi(x)\T\bigr)=u(t_0)-\tilde u\bigl(\phi(x)\T\bigr)
  \approx
  1-\T^\kappa \tilde u\circ \phi(x) \, .
$$
Thus, given \XAsympt\ and that $\psi(x)=o(x)$,
$$
  X=Ru(T)
  \approx x+\psi(x)\R - \T^\kappa x\tilde u\circ\phi(x) \, .
$$
We see that for $\R$ and $\T$ to contribute to $X$ (that is, to 
find the limiting behavior of $R$ and $T$ conditioned on $X>x$),
we should have $\psi(x)$ and $x\tilde u\circ\phi(x)$ of the same order
of magnitude (otherwise, one of the terms would dominate the other one, and
either $\R$ or $\T$ would be lost in the asymptotic). Therefore, we should
define $\phi$ by requiring $\tilde u\circ\phi(x)\sim \psi(x)/x$ as $x$
tends to infinity. We would then obtain
$$
  X\approx x+\psi(x) (\R-\T^\kappa)
$$
and the condition that $X>x$ translates into $\R>\T^\kappa$. It remains us to
formalize this sketch and turn it into a proof.

As most of the time with asymptotic analysis of integrals involving
regularly varying functions,
we will need a little more than just the definition, namely Potter's bounds.
To say that $\tilde u$ is regularly varying at $0$ with positive index $\kappa$
means that $u(1/t)$ is regularly varying with index $-\kappa$ at infinity. 
Potter's bounds are that $\tilde u(1/t)/\tilde u(1/s)$ is sandwiched between
quantities of the form $A^{\pm 1}(s/t)^{\kappa\pm\eta}$ where the real 
number $A$ can be chosen as close to $1$ as one wants, $\eta$ is positive 
and we take is less than $\kappa$, and the sandwich is good whenever $t$ 
and $s$ are large enough (see Bingham, Goldie and Teugels, 1989, \S 1.5). 
Consequently, given an $A$ greater than $1$, and a positive $\eta$, the ratio
$\tilde u(s)/\tilde u(t)$ is sandwiched between quantities of the form 
$A^{\pm 1}(s/t)^{\kappa\pm\eta}$ whenever $s$ and $t$ are small enough ---~say less
than some $\epsilon_0$.

The proof of Theorem \mainTh\ has two steps, tightness and convergence,
which are disguised as asymptotic analysis of some integrals.

We will use repeatedly that, since $u(t_0)=1$,
$$
  u(t)=1-\tilde u(t-t_0) \, .
$$

\medskip

\noindent {\it Step 1. Convergence.} Let $f$ be a nonnegative continuous 
function on $\RR^2$, 
whose support is a compact subset of $(\RR\setminus\{\, 0\,\})^2$. Consider
the integral
$$
  I(x)=\int f\Bigl({r-x\over \psi(x)},{t-t_0\over \phi(x)}\Bigr) 
  \One\{\, ru(t)>x\,;\, t>t_0\,\}
  g(t) \d H(r)\d t \, .
$$
This integral is 
$$
  {\rm E}\biggl( f\Bigl({R-x\over \psi(x)},{T-t_0\over \phi(x)}\Bigr) 
  \One\{\, X>x\,;\, T>t_0\,\}\biggr) \, ,
$$
that is, the conditional expectation given $X>x$ and $T>t_0$ multiplied 
by $P\{\, X>x\,;\, T>t_0 \,\}$.
The change of variables consisting in substituting $r$ for $(r-x)/\psi(x)$ 
and $t$ for $(t-t_0)/\phi(x)$ yields
$$\displaylines{\qquad
  I(x)=
  \int f(r,t)\One\bigl\{\, \bigl(x+r\psi(x)\bigr)
  u\bigl(t_0+t\phi(x)\bigr) >x\,;\, t>0 \,\bigr\}
  \hfill\cr\hfill
  \tilde g\bigl(t\phi(x)\bigr) \d H\bigl(x+r\psi(x)\bigr) 
  \d t\, .
  \qquad\equa{ProofA}\cr}
$$
Since $f$ has compact support which excludes the $0$-coordinates, this integral
is in fact an integral over a compact subset of $\RR^2$ which excludes
$r=0$ and $t=0$. Since $r$ and $t$ are now in a 
compact set which excludes $0$, the regular variation properties of the 
various functions yield
$$
  u\bigl(t_0+t\phi(x)\bigr)
  = 1-\tilde u\bigl(t\phi(x)\bigr)
  = 1-t^\kappa \tilde u\circ \phi(x)\bigl(1+o(1)\bigr)
$$
and
$$
  \tilde g\bigl(t\phi(x)\bigr)=t^\tau \tilde g\circ\phi(x)\bigl(1+o(1)\bigr)
$$
as $x$ tends to infinity, and both  $o(1)$ are uniform in $t$ such that 
$(r,t)$ is in the support of $f$ ---~again, because we excluded the axis 
of $\RR^2$. Thus, since $\psi(x)=o(x)$, we have
$$\eqalign{
 \bigl(x+r\psi(x)\bigr) u\bigl(t_0+t\phi(x)\bigr)
  &{}=\bigl(x+r\psi(x)\bigr) \Bigl(1-\tilde u\bigl(t\phi(x)\bigr)\Bigr) \cr
  &{}=x+r\psi(x)-x\bigl(1+o(1)\bigr) \tilde u\bigl(t\phi(x)\bigr) \cr
  &{}=x+r\psi(x)-xt^\kappa \tilde u\circ \phi(x)\bigl(1+o(1)\bigr) \, .\cr
  }
$$
Thus, referring to part of the integrand in \ProofA, and using the definition 
of $\phi$,
$$\displaylines{\qquad
  \One\bigl\{\, \bigl(x+r\psi(x)\bigr)
  u\bigl(t_0+t\phi(x)\bigr) >x\,;\, t>0\,\bigr\}
  \hfill\equa{ProofB}\cr\noalign{\vskip 5pt}\hfill
  \eqalign{
    {}={}&\One\bigl\{\, r\psi(x)-t^\kappa x\tilde u\circ\phi(x)
          \bigl(1+o(1)\bigr) >0\,;\, t>0\,\bigr\} \cr
    {}={}&\One\Bigl\{\, \psi(x)\Bigl(r-t^\kappa\bigl(1+o(1)\bigr)\Bigr)>0
                         \,;\, t>0 \,\Bigr\} \cr
    {}={}&\One\bigl\{\, r>t^\kappa\bigl(1+o(1)\bigr)\,;\, t>0\,\bigr\} \, .
  }\qquad\cr}
$$
If $x$ is large, 
the previous display shows that the indicator function in \ProofB\ can be 
sandwiched between functions
$$
  \One\{\, r>(1-\epsilon)t^\kappa\,;\, t>0\,\} \, ,
$$
(take $\epsilon$ positive for an upper bound, $\epsilon$ negative for an 
lower bound). That allows us to sandwich $I(x)$ between integrals of the form
$$\displaylines{\quad
  I_\epsilon(x)
  =\int f(r,t)\One\{\, r>(1-\epsilon)t^\kappa\,;\, t>0\,\}
  \tilde g\circ \phi(x)\phi(x)
  \hfill\cr\hfill
  t^\tau \d t \d H\bigl(x+\psi(x)r\bigr) \, ,\qquad
  \equa{ProofC}\cr}
$$
provided $x$ is large enough; thus for $\epsilon$ positive and $x$ large
enough,
$$
  (1-\epsilon)I_{-\epsilon}(x)\leq I(x)\leq (1+\epsilon)I_\epsilon(x) \, .
$$
The measure $\d H(x+\psi(x)r)/\overline H(x)$ converges vaguely to a measure
with density $e^{-r}$ with respect to the Lebesgue measure ---~note that we are
using vague convergence of measure, so that $r$ has to remain in 
a compact set, which is why we took $f$ having a compact support with 
respect to both variables $r$ and $t$. Consequently, we obtain
$$\displaylines{\qquad
  \lim_{x\to\infty} 
  {I_\epsilon(x)\over\phi(x)\tilde g\circ\phi(x)\overline H(x)}
  \hfill\cr\hfill
  {}=\int f(r,t)\One\{\, r>(1-\epsilon)t^\kappa\,;\, t>0\,\} 
  t^\tau \d t \, e^{-r}\d r 
  \qquad\equa{ProofD}\cr}
$$
as $x$ tends to infinity. Since $\epsilon$ is arbitrary, combining \ProofC\ and
\ProofD\ yield
$$
  \lim_{x\to\infty} {I(x)\over\phi(x)\tilde g\circ\phi(x)\overline H(x)}
  =
  \int f(r,t)\One\{\, r>t^\kappa\,;\, t>0\,\} t^\tau \d t\, e^{-r}\d r \, .
$$

\medskip

\noindent{\it Step 1+1/2. Refinement.} In step 1, the function $f$ is
supported in $(\RR\setminus\{\,0\,\})^2$. To prove vague convergence
of the distribution as distribution on $\RR^2$, we need to allow for
compact support in the entire $\RR^2$, not excluding the axes. To make
this extension, it suffices to show that there is no mass accumulation
along the axes $\{\,0\,\}\times\RR$ and $\RR\times\{\, 0\,\}$. Thus, setting
$$
  J_{1,\epsilon}(x)
  = \probb{ {|R-x|\over\psi(x)}\leq\epsilon \,;\, Ru(T)>x\,;\, T>t_0}
$$
and
$$
  J_{2,\epsilon}(x)
  = \probb{ {|T-t_0|\over \phi(x)} \leq\epsilon \,;\, Ru(T)>x\,;\, T>t_0} \, ,
$$
we need to prove that for $j=1,2$,
$$
  \lim_{\epsilon\to 0}\limsup_{x\to\infty} {J_{j,\epsilon}\over
  \phi(x)\tilde g\circ\phi(x)\overline H(x)} = 0 \, .
   \eqno{\equa{ProofE}}
$$
To do this, we have, for $x$ large enough, that $J_{1,\epsilon}(x)$ is at most
\hfuzz=3pt
$$\displaylines{
  \probb{ |R-x|\leq\epsilon\psi(x) \,; \, 
    u(T) >{x\over x+\epsilon\psi(x)}\,;\, T>t_0 } 
  \hfill\cr\quad
  {}\leq \Bigl( \oH\bigl(x-\epsilon\psi(x)\bigr)
    -\oH\bigl(x+\epsilon\psi(x)\bigr)\Bigr)
  %\hfill\cr\qquad\ 
    \probb{ u(T)>1-2\epsilon {\psi(x)\over x} \, ;\, T>t_0} 
  \hfill\cr\quad
  {}\leq\overline H(x) (e^\epsilon-e^{-\epsilon}) \bigl(1+o(1)\bigr)
    \probb{ \tilde u(T-t_0)<2\epsilon {\psi(x)\over x}\,;\, T>t_0} \, \cr}
$$
\hfuzz=0pt
the last inequality coming from $\overline H\in\Gamma(\psi)$, 
the definition of $\tilde u$ and that $u(t_0)=1$. 
But since $\tilde u$ is regularly varying with index $\kappa$,
$$\eqalign{
  \tilde u\bigl( (2\epsilon)^{1/\kappa}\phi(x)\bigr)
  &{}\sim 2\epsilon \tilde u\circ \phi(x)\cr
  &\sim 2 \epsilon\psi(x)/x\cr}
$$
as $x$ tends to infinity. Consequently, for $x$ large enough,
$$\displaylines{\qquad
   \probb{\tilde u(T-t_0)<2\epsilon{\psi(x)\over x}\,;\, T>t_0}
  \hfill\cr\noalign{\vskip 3pt}\hfill
  \eqalign{
    {}\leq{}& \Prob{ \tilde u(T-t_0)<\tilde u\bigl( (4\epsilon)^{1/\kappa}\phi(x)
              \bigr)\,;\, T>t_0} \cr
    {}\leq{}& \prob{|T-t_0|<(8\epsilon)^{1/\kappa}\phi(x)\,;\, T>t_0} \, ,\cr}
  \qquad\cr}
$$
the last inequality coming from the fact that a regularly varying function
of positive index is asymptotically equivalent to a monotone function
---~see Bingham, Goldie and Teugels (1989, \S 1.5.2).\note{first part of assumption 2. check}

Note that for any $\theta$ positive,
$$
  \prob{ 0<T-t_0\leq \theta\phi(x)}
  \sim \phi(x)\tilde g\circ\phi(x)\int_0^\theta y^{-\tau} \d y \, ,
  \eqno{\equa{ProofF}}
$$
as $x$ tends to infinity, because this probability is
$$
  \int_0^{\theta\phi(x)} \tilde g(s)\d s
  = \phi(x)\tilde g\circ\phi(x) \int_0^\theta  
  {\tilde g\bigl(s\phi(x)\bigr)\over \tilde g\circ\phi(x)} \d s
$$
and $\tilde g$ is regularly varying with index $\tau>-1$. Thus, combining the
various bounds, we have, for $x$ large enough,
$$
  J_{1,\epsilon}(x)
  \leq 2 \oH (x)(e^\epsilon-e^{-\epsilon}) \phi(x)\tilde g\circ\phi(x) 
  \int_0^{(16\epsilon)^{1/\kappa}} y^{-\tau}\d y
$$
and this proves \ProofE\ for $j=1$. 

To prove \ProofE\ for $j=2$, we see that for $x$ large enough,
$$\eqalign{
  J_{2,\epsilon}(x)
  &{}\leq \prob{ 0<T-t_0\leq\epsilon\phi(x) \,;\,  R>x} \cr
  &{}=\prob{0<T-t_0\leq \epsilon \phi(x)} \oH(x)\, . \cr
  }
$$
Then, we use \ProofF\ to bound $J_{2,\epsilon}(x)$,
establishing \ProofE\ for $j=2$.

Combined with Step 1, this shows that for any nonnegative continuous 
compactly supported function $f$ on $\RR^2$
$$\displaylines{\quad
  \E \biggl( f\Bigl({R-x\over\psi(x)},{T-t_0\over \phi(x)}\Bigr)
  \One\{\, X>x\,;\, T>t_0\,\}\biggr)
  \hfill\cr\hfill
  \sim \phi(x)\tilde g\circ\phi(x) \oH(x)
  \int f(r,t) \One\{\, r>t^\kappa\,;\, t>0\,\} t^\tau \d t \, e^{-r} \d r
  \quad\cr}
$$
as $x$ tends to infinity. By writing any continuous function as the sum
of its positive and negative part, this still holds for any continuous and
compactly supported function on $\RR^2$.

\medskip

\noindent{\it Step 2. Tightness.} We now show that $(R-x)/\psi(x)$ 
and $(T-t_0)/\phi(x)$ are tight random variables under the conditional 
probability given $X>x$ and $T>t_0$.
For this purpose, given step 1 and anticipating the conclusion of the proof, 
we need to show that
$$
  \lim_{r\to\infty}\limsup_{x\to\infty}
  { \probb{ {\ds|R-x|\over\ds \psi(x)}>r \, ;\, Ru(T)>x\,;\, T>t_0} 
    \over
    \phi(x)\tilde g\circ\phi(x)\overline H(x) }
  = 0\, ,
$$
and
$$
  \lim_{t\to\infty} \limsup_{x\to\infty}
  { \probb{ {\ds |T-t_0|\over\ds \phi(x)}>t \, ;\,  Ru(T)>x\,;\, T>t_0}
    \over 
    \phi(x)\tilde g\circ\phi(x)\overline H(x) }
  = 0\, ,
  \eqno{\equa{ProofFa}}
$$
This is a bit painful, because of the absolute values involved. We will
examine the three cases obtained when `removing' the absolute values.

\noindent{\it Case 1.} Let $r$ be positive and let us bound
$$
  P_{1,r}(x) = \prob{R>x+r\psi(x) \,;\, Ru(T)>x\,;\, T>t_0} \, .
$$
For $x$ large enough, this is at most
$$\displaylines{\qquad
  \oH\bigl( x+r\psi(x)\bigr)\probb{ u(T)>{x\over x+r\psi(x)}\,;\, T>t_0 }
  \hfill\cr\hfill
  \sim \oH(x)e^{-r} 
    \probb{\tilde u(T-t_0)< r{\psi(x)\over x}\bigl(1+o(1)\bigr)\,;\, T>t_0}
    \, .
  \qquad\cr}
$$
As in step 1+1/2, using \ProofF, this is of order at most
$$
  \oH(x) \phi(x)g\circ\phi(x) e^{-r} 
  \int_0^{(4r)^{1/\kappa}} y^{-\kappa}\d y \, .
$$
Thus,
$$
  \lim_{r\to\infty} \limsup_{x\to\infty} 
  {P_{1,r}(x)\over \phi(x)\tilde g\circ\phi(x) \oH(x)} = 0 \, .
$$

\noindent{\it Case 2.} For $r$ positive, define
$$
  P_{2,r}(x)= \prob{ R<x-r\psi(x) \,;\, Ru(T)>x\,;\, T>t_0 } \, .
$$
When $x$ is large enough, $\psi(x)$ is well defined and positive. In that
range of $x$, since $|u|\leq 1$, we cannot have $Ru(T)>x$ while 
having $R<x-r\psi(x)$. Thus $P_{2,r}(x)=0$ whenever $x$ is large enough.

\noindent{\it Case 3.} The probability involved in the numerator of \ProofFa\
is
$$
  P_{3,t}(x) =\prob{ T>t_0+t\phi(x) \,;\, Ru(T)>x\,;\, T>t_0} \, .
$$
We see that
$$
  P_{3,t}(x)
  = \int_t^\infty\overline H\Bigl( {x\over u\bigl(t_0+\phi(x)s\bigr)}\Bigr)
  \phi(x) \tilde g\bigl(s\phi(x)\bigr) \d s \, .
  \eqno{\equa{ProofG}}
$$
We write 
$$
  u\bigl(t_0+s\phi(x)\bigr) = 1-\tilde u\bigl(s\phi(x)\bigr)
$$
and we use the usual arguments to bound the integral: Potter's bound
whenever we can, and ad hoc argument elsewhere. This goes as follows.
We may assume that $t$ is greater than $1$.
Let $\eta$ be a (small) positive real number.
Let $\epsilon$ be small enough so that Potter's bounds
$$
  \tilde u\bigl(s\phi(x)\bigr)\geq {1\over 2}\tilde u\circ\phi(x) s^{\kappa-\eta}
$$
and
$$
  \tilde g\bigl(s\phi(x)\bigr) \leq 2 \tilde g\circ\phi(x) s^{\tau+\eta}
  \eqno{\equa{ProofH}}
$$
apply on the range $1\leq t\leq s\leq \epsilon/\phi(x)$. We then have, on that
range of $s$ (provided $\epsilon$ was chosen small enough),
$$\eqalignno{
  {1\over u\bigl(t_0+s\phi(x)\bigr)} 
  = { 1\over 1-\tilde u\bigl(s\phi(x)\bigr) }
  &{}\geq 1+{1\over 4}\tilde u\bigl(s\phi(x)\bigr)\cr
  &{}\geq 1+{1\over 8}\tilde u\circ \phi(x)s^{\kappa-\eta} \, .&\equa{ProofK}\cr
  }
$$ 
Referring to part of the integral \ProofG, using the definition 
of $\phi$, \ProofH\ and \ProofK, we have then
$$\displaylines{
  \int_t^{\epsilon/\phi(x)}\overline H\Bigl( {x\over u\bigl(t_0+s\phi(x)\bigr)}
  \Bigr) \tilde g\bigl(s\phi(x)\bigr) \phi(x) \d s 
  \hfill\cr\hfill
  {}\leq 2\phi(x)\tilde g\circ\phi(x) \int_t^{\epsilon/\phi(x)}
               \overline H\Bigl( x+{1\over 16} \psi(x)s^{\kappa-\eta}\Bigr) 
               s^{\tau+\eta} \d s\cr
 }
$$
Using the first statement of Lemma 5.1 in Foug\`eres and Soulier (2010) 
(note we can take $C=2$ in that Lemma, which we do here), this upper
bound is at most
$$
  4\phi(x)\tilde g\circ \phi(x)\overline H(x)\int_t^\infty 
  {s^{\tau+\eta}\over \bigl(1+(s^{\kappa-\eta}/16)\bigr)^p} \d s \, ,
  \eqno{\equa{ProofL}}
$$
where $p$ is taken large enough so that the integral converges.

We now work on the easy part of the integral \ProofG, namely, that for $s$
between $\epsilon/\phi(x)$ and $\infty$. Given how this integral was obtained,
this part corresponds to $T> t_0+\epsilon$, and it is at most (again, 
provided we choose $\epsilon$ small enough)
$$
  \prob{R u(t_0+\epsilon/2)>x} 
  = \overline H\Bigl({x\over u(t_0+\epsilon/2)}\Bigr) \, .
  \eqno{\equa{ProofM}}
$$
We now claim that if $c>1$ (think of $c$ as $1/u(t_0+\epsilon/2)$), then
$$
  \overline H(cx)=o\bigl(\overline H(x)\phi(x)\tilde g\circ\phi(x)\bigr) \, .
  \eqno{\equa{ProofN}}
$$
Indeed, using the second statement of Lemma 5.1 in Foug\`eres and 
Soulier (2010), for any positive $p$ we have
$$
  \overline H(cx)
  \leq \Bigl({\psi(x)\over x}\Bigr)^p\overline H(x)
$$
provided $x$ is large enough (note that we can take $C=1$ in their inequality:
it suffices to divide their $p$ by $2$ and see that their $C$ 
times $(\psi(x)/x)^{p/2}$ tends to $0$ and is less than $1$ for $x$ large 
enough). Thus, to prove \ProofN, we have to show that for any $p$ large enough
$$
  \Bigl({\psi(x)\over x}\Bigr)^p
  = o\bigl(\phi(x)\tilde g\circ\phi(x)\bigr) \, .
$$
But this comes from viewing $\phi(x)\tilde g\circ\phi(x)$ has a function
of $\psi(x)/x$ which is then regularly varying of index $(\tau+1)/\kappa$ in
that argument.

Now, combining \ProofM\ and \ProofN, we obtain that, referring to part 
of \ProofG\
$$
  \int_{\epsilon/\phi(x)}^\infty\overline 
  H\Bigl( {x\over u\bigl(t_0+\phi(x)s\bigr)}\Bigr)
  \tilde g\bigl(\phi(x)s\bigr) \phi(x) \d s 
  = o\bigl(\phi(x)\tilde g\circ\phi(x)\overline H(x)\bigr)
$$
as $x$ tends to infinity. Combined with \ProofL, and referring to \ProofG\
this shows that
$$
  \limsup_{x\to\infty}
    { P_{3,t}(x)\over\phi(x)\tilde g\circ\phi(x)\overline H(x)}
  \leq 4 \int_t^\infty {s^{\tau+\eta}\over \bigl(1+(s^{\kappa-\eta}/16)\bigr)^p}\d s
  \, ,
$$
and, therefore,
$$
  \lim_{t\to\infty}\limsup_{x\to\infty}
  { P_{3,t}(x)\over\phi(x)\tilde g\circ\phi(x)\overline H(x)}
  =0 \, .
$$

To conclude the proof, combining steps 1, 1+1/2 and 2, we obtain that
$$\displaylines{\qquad
  \prob{X>x\,;\, T>t_0}
  \hfill\cr\noalign{\vskip 5pt}\hfill
  \eqalign{
  {}\sim{}& \phi(x)\tilde g\circ\phi(x) \oH(x) 
            \int\One\{\, r>t^\kappa\,;\, t>0\,\}
            t^{\tau} e^{-r} \d t\d r \cr
  {}\sim{}&  \phi(x)\tilde g\circ\phi(x) \oH(x) 
    {1\over\kappa}\Gamma\Bigl({1+\tau\over\kappa}\Bigr) \cr}
  \qquad\cr}
$$
as $x$ tends to infinity. Then, step 2 implies that the conditional 
distribution of 
$$
  \Bigl( {R-x\over\psi(x)},{T-t_0\over\phi(x)}\Bigr)
$$
given $X>x$ and $T>t_0$ is tight, and step 1 proves that it converges to the 
limit given in Theorem \mainTh.\hfill\qed

\bigskip

\noindent{\bf Acknowledgements.} The authors thank Anne-Laure Foug\`eres
for making their collaboration possible.

%%%%%%%%%%%%%%%%%%%%%%%%%%%%%%%%%%%%%%%%%%%%%%%%%%%%%%%%%%%%%%%%%%%%%%
% References
%%%%%%%%%%%%%%%%%%%%%%%%%%%%%%%%%%%%%%%%%%%%%%%%%%%%%%%%%%%%%%%%%%%%%%

\bigskip

\noindent{\bf References.}
\medskip

{\leftskip=\parindent \parindent=-\parindent
 \par

S.M.\ Berman (1983). Sojourns and extremes of Fourier sums and series with 
random coefficients, {\sl Stochastic Process.\ Appl.}, 15, 213--238. 

N.H.\ Bingham, C.M.\ Goldie, J.L.\ Teugels (1989). {\sl Regular Variation},
2nd ed. Cambridge University Press.

B.\ Das, S.I.\ Resnick (2011). Detecting a conditional extreme value model,
{\sl Extremes}, 14, 29--61.

B.\ Das, S.I.\ Resnick (2011). Conditioning on an extreme component: Model
consistency and regular variation on cones, Bernoulli, 17, 226--252.

L.\ de Haan (1970). {\sl On Regular Variation and its Application to the Weak
Convergence of the Sample Extremes}, Mathematisch Centrum, Amsterdam.

R.M.\ Dudley (1989). {\sl Real Analysis and Probability}, Chapman\& Hall.

A.-L.\ Foug\`eres, Ph.\ Soulier (2010). Limit conditional distributions for 
bivariate vectors with polar representation, {\sl Stoch.\ Models}, 26, 54-77.

A.-L.\ Foug\`eres, Ph.\ Soulier (2012) Estimation of conditional laws given 
an extreme component, {\sl Extremes}, 15, 1--34.

E.\ Hashorva (2012). Exact tail asymptotics in bivariate scale mixture 
models, {\sl Extremes}, 15, 109--128.

J.E.\ Heffernan, S.I.\ Resnick (2007). Limit laws for random vectors with an 
extreme component, {\sl Ann.\ Appl.\ Probab.}, 17, 537--571.

J.E.\ Heffernan, J.A.\ Tawn (2004). A conditional approach for multivariate 
extreme values, {\sl J.\ R.\ Stat.\ Soc.\ Ser.\ B Stat.\ Methodol.}, 66, 
497--546.

M.M.\ Merschaert, H.-P.\ Scheffler (2001). {\sl Limit Distributions for Sums
of Independent Random Vectors}, Wiley.

S.I.\ Resnick (1987). {\sl Extreme Values, Regular Variation, and Point
Processes}, Springer.

M.I.\ Seifert (2012). On conditional extreme values of random vectors with
polar representation, {\sl preprint}.

V.S.\ Vladimirov, Yu.N.\ Drosinov, B.I.\ Zavialov (1988). {\sl Tauberian
Theorems for Generalized Fuctions}, Kluwer.

}

\bigskip\bigskip

\setbox1=\vbox{\halign{#\hfil &\hskip 40pt #\hfill\cr
  Ph.\ Barbe                  & M.I.\ Seifert\cr
  90 rue de Vaugirard         & Helmut-Schmidt Universit\"at\cr
  75006 PARIS                 & Holstenhofweg 85\cr
  FRANCE                      & 22043 Hamburg \cr
  philippe.barbe@math.cnrs.fr & Germany \cr
                              & miriam.seifert@hsu-hh.de \cr}}
\box1

\bigskip

\bye